\documentclass{article}
\usepackage{spconf,amsmath,graphicx}
\usepackage{amsmath,amssymb,amsthm}%
\begin{document}
\title{Compressive Joint Angular-Frequency Power Spectrum Estimation}
\vspace{-1mm}
\name{Dyonisius Dony Ariananda and Geert Leus}
\vspace{-1mm}
\address{Faculty of EEMCS, Delft University of Technology, The Netherlands\\
\{d.a.dyonisius, g.j.t.leus\}@tudelft.nl.}
\maketitle
\thispagestyle{empty}
\vspace{-1mm}
\begin{abstract}
We introduce a new compressive power spectrum estimation approach in both frequency and direction of arrival (DOA). Wide-sense stationary signals produced by multiple uncorrelated sources are compressed in both the time and spatial domain where the latter compression is implemented by activating only some of the antennas in the underlying uniform linear array (ULA).
We sample the received signal at every active antenna at sub-Nyquist rate, 
compute both the temporal and spatial correlation functions between the sub-Nyquist rate samples, 
and apply least squares to 
reconstruct the full-blown two-dimensional power spectrum matrix where the rows and columns 
correspond to the frequencies and the angles, respectively. 
This is possible under the full column rank condition of the system matrices and without applying any sparsity constraint on the signal statistics. Further, 
we can estimate the DOAs of the sources by locating the peaks of the angular power spectrum. 
We can theoretically estimate the frequency bands and the DOAs of more uncorrelated sources than active sensors using sub-Nyquist sampling.
\end{abstract}
\vspace{-2mm}
\section{Introduction}\label{introduction}
\vspace{-1mm}
Compressive sampling and multi-coset sampling 
have 
drawn a lot of interest from the signal processing community due to the possibility to reconstruct a signal sampled at sub-Nyquist rate with no or little information loss under the constraint that the signal is sparse in a particular basis~\cite{Candes,Yonina}.
All these works on sub-Nyquist sampling are important especially when it is needed to relax the requirements on the analog-to-digital converters (ADCs). 
For a wide-sense stationary (WSS) signal,
it has also been shown that perfect reconstruction of its second-order statistics from sub-Nyquist rate samples is theoretically possible even without sparsity constraint~\cite{TSP12}. This invention is important for some applications, such as wideband spectrum sensing for cognitive radio, where only perfect reconstruction of the temporal auto-correlation function 
is required instead of the signal itself.
The principle of reconstructing the temporal auto-correlation 
function of a signal from the time-domain compressive measurements has in a dual form also been proposed in the spatial domain. Given a linear antenna array,~\cite{Nested} and~\cite{Coprime} show that if the locations of the antennas are arranged according to a nested 
or coprime array, the spatial correlation values between the outputs of the antennas in the array can be used to generate the spatial correlation values between the outputs of the antennas in the virtual array or difference co-array (which is uniform in this case) which generally has more antennas and a larger aperture than the actual array. 
This enhances the degrees of freedom 
and allows~\cite{Nested} and~\cite{Coprime} to estimate the direction of arrival (DOA) of more uncorrelated sources than sensors. 
The 
minimum redundancy array (MRA) of~\cite{Moffet} can also be used to produce this feature 
but in a more optimal way. This has been exploited by~\cite{Siavash} to perform compressive angular power spectrum reconstruction. The advantage offered by the nested and coprime arrays over the MRA however, is the possibility to derive a closed-form expression for the array geometry and the achievable number of correlation values in the resulting uniform difference co-array. 
In the aforementioned concept, the spatial compression is performed in the sense that we select a subset of antennas from a uniform linear array (ULA).

In this paper, we jointly reconstruct both the frequency-domain and angular-domain power spectrum using compressive samples. 
We use a ULA as the underlying array and activate only some of its antennas
leading to a spatial-domain compression. The received signal at each active antenna is then sampled at sub-Nyquist-rate using multi-coset sampling. 
Next, we compute all the correlation values between the resulting sub-Nyquist rate samples at all active antennas both in the time domain and the spatial domain and use them 
to reconstruct the 
two-dimensional (2D) power spectrum matrix where each row gives the power spectrum in the frequency domain for a given angle and where each column contains the power spectrum 
in the angular domain for a given frequency. Further, 
we can estimate the DOA of the sources active at each frequency by locating the peaks in the angular power spectrum. This 2D power spectrum reconstruction can be done for more uncorrelated sources than active sensors without any sparsity constraint on the true power spectrum. 
\vspace{-2mm}
\section{Preliminaries}\label{preliminary}
\vspace{-1mm}
First, consider a ULA having $N_s$ antennas receiving signals from $K$ uncorrelated WSS sources. We assume that the distance between the sources and the ULA is large enough compared to the length of the ULA and thus the wave incident on the ULA is assumed to be planar and the sources can be assumed as point sources. We also assume that the inverse of the bandwidth of the aggregated incoming signals is larger than the propagation delay across the ULA, which allows us to represent the delay between the antennas 
as a phase shift. Based on these assumptions, we can write the ULA output as
\vspace{-2mm}
\begin{equation}
{\bf x}(t)=\sum_{q=1}^Q{\bf a}(\theta_q)s_q(t)+{\bf n}(t)={\bf A}{\bf s}(t)+{\bf n}(t)
\vspace{-2mm}
\label{eq:ULA_output}
\end{equation}
where ${\bf x}(t)$ is the $N_s \times 1$ output vector containing the received signal at the $N_s$ antennas of the ULA, ${\bf n}(t)$ is the $N_s \times 1$ additive white Gaussian noise vector, ${\bf s}(t)=[s_1(t),s_2(t),\dots,s_Q(t)]^T$ is the $Q \times 1$ extended source vector with $s_q(t)$ the incoming signal from the investigated angle 
$\theta_q$, and ${\bf A}=\left[{\bf a}(\theta_1), {\bf a}(\theta_2), \dots, {\bf a}(\theta_{Q})\right]$ is the $N_s \times Q$ extended array manifold matrix with ${\bf a}(\theta_q)$ the $N_s \times 1$ array response vector containing the phase shifts experienced by $s_q(t)$ at each element of the ULA. Note that $\{\theta_q\}_{q=1}^Q$ is known and might only approximately contain the actual DOAs of the $K$ sources. 
We generally assume that ${\bf n}(t)$ and ${\bf s}(t)$ are uncorrelated, that the impact of the wireless channel has been taken into account in ${\bf s}(t)$, and that the noises at different antennas are uncorrelated with variance $\sigma_n^2$, i.e., $E[{\bf n}(t){\bf n}^H(t)]=\sigma_n^2 {\bf I}_{N_s}$, with ${\bf I}_{N_s}$ the $N_s \times N_s$ identity matrix. 
We consider the first element of the ULA as a reference point and express the array response vector ${\bf a}(\theta_q)$ as ${\bf a}(\theta_q)=[1, a(\theta_q)^d, a(\theta_q)^{2d}, \dots, a(\theta_q)^{(N_s-1)d}]^T$, where 
$a(\theta_q)=\text{exp}\left(j 2\pi\text{sin}(\theta_q)\right)$ 
and $d$ is the distance between two consecutive antennas in wavelengths, which is set to $d \leq 0.5$ in order to prevent spatial aliasing. 

In order to simplify the further analysis, we introduce ${\bf x}[m]={\bf x}(mT)$, ${\bf n}[m]={\bf n}(mT)$, and ${\bf s}[m]={\bf s}(mT)$ as a digital representation of ${\bf x}(t)$, ${\bf n}(t)$, and ${\bf s}(t)$, respectively, where $1/T$ is the Nyquist sampling rate at every ADC associated with each antenna. We then 
collect the output vectors ${\bf x}[m]$ at $N_t$ consecutive sample indices into 
the $N_s \times N_t$ matrix ${\bf X}[n]$, for $n=0,1,\dots,N_{n-1}$, as 
${\bf X}[n]=\left[{\bf x}[nN_t],{\bf x}[nN_t+1],\dots,{\bf x}[(n+1)N_t-1]\right]$ 
and write ${\bf X}[n]$ as
\begin{equation}
{\bf X}[n]={\bf AS}[n]+{\bf N}[n]
\label{eq:X_discrete}
\end{equation}
where ${\bf N}[n]$ is similarly defined as ${\bf X}[n]$ and 
the $Q \times N_t$ matrix ${\bf S}[n]$ is given by ${\bf S}[n]=[{\bf s}[nN_t],{\bf s}[nN_t+1],\dots,{\bf s}[(n+1)N_t-1]]$. Let us also write the $Q\times 1$ vector ${\bf s}[nN_t+i]$ as ${\bf s}[nN_t+i]=[{s}_1[nN_t+i],{s}_2[nN_t+i],\dots,{s}_Q[nN_t+i]]^T$ with ${s}_q[m]={s}_q(mT)$ a digital representation of $s_q(t)$. 
\vspace{-1mm} 
\section{Time-Domain and Spatial-Domain Compression}\label{time-spatial}
\vspace{-1mm}
In this section, we introduce the compression operations on the output matrix ${\bf X}[n]$ both in the spatial domain and time domain. The spatial-domain compression is implemented by activating only $M_s$ out of $N_s$ available antennas in the ULA leading to a possibly non-ULA of less active antennas than sources. Further, in the receiver branches associated with the $M_s$ active antennas, time-domain compression is performed by sampling the received analog signal 
at sub-Nyquist-rate 
using the multi-coset sampling principle discussed in~\cite{Yonina}, 
which can be implemented using the practical sampling device proposed in~\cite{TSP12}. Here, the multi-coset sampling process is represented by selecting only $M_t$ out of $N_t$ time samples.  

We first introduce the $M_s \times N_s$ spatial-domain selection matrix ${\bf C}_s$, which is formed by selecting $M_s$ rows of ${\bf I}_{N_s}$. 
Here, the indices of the selected rows of ${\bf I}_{N_s}$ used to construct ${\bf C}_s$ correspond to the indices of the $M_s$ active antennas selected from the $N_s$ available antennas in the ULA. Based on~\eqref{eq:X_discrete}, the $N_s \times N_t$ matrix ${\bf X}[n]$ is then compressed in the spatial-domain by ${\bf C}_s$ 
leading to the $M_s \times N_t$ matrix 
\vspace{-1mm}
\begin{equation}
{\bf Y}[n]={\bf C}_s{\bf X}[n]
\stackrel{\Delta}{=}{\bf B}{\bf S}[n]+{\bf M}[n]
\vspace{-1mm}
\label{eq:Y_as_X}
\end{equation}
where ${\bf Y}[n]=[{\bf y}[nN_t],{\bf y}[nN_t+1],\dots,{\bf y}[(n+1)N_t-1]]$ with ${\bf y}[nN_t+l]=[y_1[nN_t+l],y_2[nN_t+l],\dots,y_{M_s}[nN_t+l]]^T$, ${\bf B}=\left[{\bf b}(\theta_1), {\bf b}(\theta_2), \dots, {\bf b}(\theta_{Q})\right]$ is the $M_s \times Q$ array response matrix with ${\bf b}(\theta_q)={\bf C}_s{\bf a}(\theta_q)$ the $M_s \times 1$ array response vector associated with the $M_s$ activated antennas, the $M_s \times N_t$ matrix ${\bf M}[n]$ is given by ${\bf M}[n]=\left[{\bf m}[nN_t],{\bf m}[nN_t+1],\dots,{\bf m}[(n+1)N_t-1]\right]$, and ${\bf m}[m]$ is the $M_s \times 1$ discrete noise vector given by ${\bf m}[m]={\bf C}_s{\bf n}[m]$. Observe that ${\bf m}[m]$ generally has correlation matrix $E\left[{\bf m}[m]{\bf m}^H[m']\right]=\sigma_n^2{\bf I}_{M_s}\delta[m-m']$.
The next step is to introduce the $M_t \times N_t$ time-domain selection matrix ${\bf C}_t$ formed by selecting $M_t$ rows of the $N_t \times N_t$ identity matrix ${\bf I}_{N_t}$, and further compress ${\bf Y}[n]$ in~\eqref{eq:Y_as_X} in the time domain, leading to the $M_s \times M_t$ matrix
\begin{equation}
{\bf Z}[n]={\bf Y}[n]{\bf C}_t^T.
\label{eq:Yprime_as_Y}
\end{equation}
\vspace{-7mm}
\section{Power Spectrum Reconstruction}\label{psd_reconstruction}
\vspace{-1mm}
Denote the $j$-th row of ${\bf Z}[n]$ and ${\bf Y}[n]$ in~\eqref{eq:Yprime_as_Y} as ${\bf z}^T_j[n]$ and ${\bf y}^T_j[n]$, respectively, and write the $M_t \times 1$ vector ${\bf z}_j[n]$ in terms of its elements as ${\bf z}_j[n]=[z_{j,1}[n],z_{j,2}[n],\dots,z_{j,M_t}[n]]^T$ 
and the $N_t \times 1$ vector ${\bf y}_j[n]$ as ${\bf y}_j[n]=[y_j[nN_t],y_j[nN_t+1],\dots,y_j[(n+1)N_t-1]]^T$. This allows us to rewrite the time-domain compression in~\eqref{eq:Yprime_as_Y} in terms of the row vectors of ${\bf Z}[n]$ and ${\bf Y}[n]$, i.e.,
\vspace{-1mm}
\begin{equation}
{\bf z}_j[n]={\bf C}_t{\bf y}_j[n],\quad j=1,2,\dots,M_s.
\vspace{-1mm}
\label{eq:y_jprime_as_y_j}
\end{equation}
Using~\eqref{eq:y_jprime_as_y_j}, our next step is to calculate the correlation matrix between ${\bf z}_i[n]$ and ${\bf z}_j[n]$ for all $i,j=1,2,\dots,M_s$ as
\vspace{-1mm}
\begin{equation}
{\bf R}_{z_i,z_j}=
E\left[{\bf C}_t{\bf y}_i[n]{\bf y}_j[n]^H{\bf C}_t^H\right]={\bf C}_t{\bf R}_{y_i,y_j}{\bf C}_t^H.
\vspace{-1mm}
\label{eq:R_yi_yj_prime}
\end{equation}
In practice, the expectation operator in~\eqref{eq:R_yi_yj_prime} can be estimated by taking an average over $N_n$ available matrices ${\bf Z}[n]$. After cascading all columns of ${\bf R}_{z_i,z_j}$ into the $M_t^2 \times 1$ vector vec$({\bf R}_{z_i,z_j})$ and by taking into account the fact that ${\bf C}_t$ is a real matrix, we can express vec$({\bf R}_{z_i,z_j})$ based on~\eqref{eq:R_yi_yj_prime} as
\begin{equation}
\text{vec}({\bf R}_{z_i,z_j})
=({\bf C}_t\otimes{\bf C}_t)\text{vec}({\bf R}_{y_i,y_j})
\label{eq:vec_R_yi_yj_prime}
\end{equation}
where vec$(.)$ is the operator that cascades all columns of a matrix into a single column vector and $\otimes$ represents the Kronecker product operation. 
Up to this stage, let us recall that $\{s_q(t)\}_{q=1}^Q$ in~\eqref{eq:ULA_output} are WSS processes since we have $K$ WSS sources. Based on this fact, as well as~\eqref{eq:Y_as_X}, it is obvious that the elements of 
${\bf y}_j[n]$ in~\eqref{eq:y_jprime_as_y_j} also form a WSS sequence. This means that the $N_t \times N_t$ matrix ${\bf R}_{y_i,y_j}$ in~\eqref{eq:vec_R_yi_yj_prime} 
has a Toeplitz structure allowing us to condense ${\bf R}_{y_i,y_j}$ into the $(2N_t-1) \times 1$ vector 
${\bf r}_{y_i,y_j}=[r_{y_i,y_j}[0],r_{y_i,y_j}[1],\dots,r_{y_i,y_j}[N_t-1],r_{y_i,y_j}[1-N_t],\dots,r_{y_i,y_j}[-1]]^T$
and write
\vspace{-1mm}
\begin{equation}
\text{vec}({\bf R}_{y_i,y_j})={\bf T}{\bf r}_{y_i,y_j}
\vspace{-1mm}
\label{eq:repetition}
\end{equation}
where ${\bf T}$ is a special $N_t^2 \times (2N_t-1)$ repetition matrix whose $i$-th row is given by the $((i-1+(N_t-2)\left\lfloor \frac{i-1}{N_t}\right\rfloor) \bmod (2N_t-1) +1)$-th row of the identity matrix ${\bf I}_{2N_t-1}$. By combining~\eqref{eq:vec_R_yi_yj_prime} and~\eqref{eq:repetition}, we obtain
\vspace{-1mm}
\begin{equation}
\text{vec}({\bf R}_{z_i,z_j})=({\bf C}_t\otimes{\bf C}_t){\bf T}{\bf r}_{y_i,y_j}={\bf R}_{c_t}{\bf r}_{y_i,y_j}
\vspace{-1mm}
\label{eq:vec_R_yi_yj_prime_as_T}
\end{equation}
where ${\bf R}_{c_t}=({\bf C}_t\otimes{\bf C}_t){\bf T}$ is an $M_t^2 \times (2N_t-1)$ matrix. 
Observe that it is possible to reconstruct ${\bf r}_{y_i,y_j}$ from $\text{vec}({\bf R}_{z_i,z_j})$ in~\eqref{eq:vec_R_yi_yj_prime_as_T}, for all $i,j=1,2,\dots,M_s$, 
using least squares (LS) if $M_t^2 \geq 2N_t-1$ and ${\bf R}_{c_t}$ has full column rank.

The next step is to figure out the relationship between $\{{\bf r}_{y_i,y_j}\}_{i,j=1}^{M_s}$ 
in~\eqref{eq:vec_R_yi_yj_prime_as_T} and the extended source matrix ${\bf S}[n]$ in~\eqref{eq:Y_as_X}. 
By taking into account the fact that every row 
of ${\bf Y}[n]$ and ${\bf S}[n]$ is a WSS sequence and the assumption that the extended source 
vector ${\bf s}[m]$ and the noise vector ${\bf m}[m]$ are uncorrelated, 
it is straightforward to find that the correlation matrix between ${\bf y}[nN_t+l]$ and ${\bf y}[nN_t+l']$ is given by
\vspace{-1mm}
\begin{equation}
{\bf R}_y[l-l']
={\bf B}{\bf R}_{{s}}[l-l']{\bf B}^H+\sigma_n^2{\bf I}_{M_s}\delta[l-l']
\vspace{-1mm}
\label{eq:Ry_only}
\end{equation}
for $l,l'=0,1,\dots,N_t-1$. 
Since the point sources are assumed to be uncorrelated, the elements of 
${\bf s}[m]$ 
are also uncorrelated and thus the $Q\times Q$ matrix ${\bf R}_{s}[l-l']$ is a diagonal matrix. 
By exploiting this fact and stacking all columns of the $M_s \times M_s$ matrix ${\bf R}_y[l-l']$ in~\eqref{eq:Ry_only} into the $M_s^2 \times 1$ vector vec$({\bf R}_y[l-l'])$, 
we obtain
\vspace{-1mm}
\begin{align}
&\text{vec}({\bf R}_y[l-l'])
=({\bf B}^*\odot{\bf B})\text{diag}({\bf R}_{s}[l-l']) \nonumber \\
&+\sigma_n^2\text{vec}({\bf I}_{M_s})\delta[l-l'], \quad l,l'=0,1,\dots,N_t-1,
\vspace{-1mm}
\label{eq:vec_Ry_only}
\end{align}
where $\odot$ represents the Khatri-Rao product operation. 
Let us now investigate the relationship between the elements of $\{{\bf r}_{y_i,y_j}\}_{i,j=1}^{M_s}$ in~\eqref{eq:vec_R_yi_yj_prime_as_T} and 
$\text{vec}({\bf R}_y[l-l'])$ in~\eqref{eq:vec_Ry_only}. 
We can find that $\{\text{vec}({\bf R}_y[l-l'])\}_{l,l'=0}^{N_t-1}$ is actually related to $\{{\bf r}_{y_i,y_j}\}_{i,j=1}^{M_s}$ as 
$\text{vec}({\bf R}_y[l-l'])=[r_{y_1,y_1}[l-l'],r_{y_2,y_1}[l-l'],\dots,r_{y_{M_s},y_{M_s}}[l-l']]^T$. 
Hence, we can use the elements of the reconstructed $\{{\bf r}_{y_i,y_j}\}_{i,j=1}^{M_s}$ in~\eqref{eq:vec_R_yi_yj_prime_as_T} to form $\{\text{vec}({\bf R}_y[l-l'])\}_{l,l'=0}^{N_t-1}$ in~\eqref{eq:vec_Ry_only} and then use them to reconstruct $\{\text{diag}({\bf R}_{s}[l-l'])\}_{l,l'=0}^{N_t-1}$ in~\eqref{eq:vec_Ry_only}, which can be performed using LS if $M_s^2 \geq Q$ and ${\bf B}^*\odot{\bf B}$ has full column rank.

If we combine the $Q\times 1$ vectors $\{\text{diag}({\bf R}_{s}[l-l'])\}_{l,l'=0}^{N_t-1}$ as
$\bar{\bf R}_{s}=[\text{diag}({\bf R}_{s}[0]),\text{diag}({\bf R}_{s}[1]),\dots,\text{diag}({\bf R}_{s}[N_t-1]),$ $\text{diag}({\bf R}_{s}[1-N_t]),\dots,\text{diag}({\bf R}_{s}[-1])]$, 
we can 
observe that the $q$-th row of $\bar{\bf R}_{s}$ actually corresponds to the temporal auto-correlation of the incoming signal from the investigated angle ${\theta}_q$, which can be written as ${\bf r}^T_{s_q}=[r_{s_q}[0],r_{s_q}[1],\dots,r_{s_q}[N_t-1],r_{s_q}[1-N_t],\dots,r_{s_q}[-1]]$. By defining ${\bf F}_{2N_t-1}$ as the $(2N_t-1) \times (2N_t-1)$ discrete Fourier transform (DFT) matrix, we can compute the power spectrum of ${s}_q[m]$ as 
${\bf p}_{s_q}={\bf F}_{2N_t-1}{\bf r}_{s_q}$, 
where ${\bf p}_{s_q}$ is the $(2N_t-1) \times 1$ power spectrum vector of the incoming signal from the investigated angle ${\theta}_q$. By combining 
$\{{\bf p}_{s_q}\}_{q=1}^Q$ 
into the $Q\times (2N_t-1)$ matrix $\bar{\bf P}_{s}=[{\bf p}_{s_1},{\bf p}_{s_2},\dots,{\bf p}_{s_Q}]^T$, we can write
\vspace{-1mm}
\begin{equation}
\bar{\bf P}_{s}=\bar{\bf R}_{s}{\bf F}_{2N_t-1}.
\vspace{-1mm}
\label{eq:P_stilde}
\end{equation}
Note that $\bar{\bf P}_{s}$ in~\eqref{eq:P_stilde} can be perceived as a 2D power spectrum matrix where every row of $\bar{\bf P}_{s}$ gives the power spectrum in the frequency-domain for a given investigated angle and every column of $\bar{\bf P}_{s}$ provides the power spectrum information in the angular domain for a given frequency.
\vspace{-1mm}
\section{Construction of the Compression Matrices}
\label{sec:Compression_Matrices}
\vspace{-1mm}
Recall that the 2D power spectrum matrix $\bar{\bf P}_{s}$ can be reconstructed from $\text{vec}({\bf R}_{z_i,z_j})$ in~\eqref{eq:vec_R_yi_yj_prime_as_T}, 
which contains the cross-correlations between the rows of the measurement matrix ${\bf Z}[n]$ in~\eqref{eq:Yprime_as_Y}, by solving~\eqref{eq:vec_R_yi_yj_prime_as_T} 
and~\eqref{eq:vec_Ry_only} 
using LS and then applying the DFT on the rows of the resulting matrix $\bar{\bf R}_{s}$. 
We 
now discuss the choice of the selection matrix ${\bf C}_t$ and the extended array response matrix ${\bf B}$ that ensure 
the uniqueness of the LS solution of~\eqref{eq:vec_R_yi_yj_prime_as_T} 
and~\eqref{eq:vec_Ry_only}, respectively. 

We first investigate the choice of ${\bf C}_t$ that results in a full column rank matrix ${\bf R}_{c_t}$. Since the rows of ${\bf C}_t$ and ${\bf T}$ in~\eqref{eq:vec_R_yi_yj_prime_as_T} are formed by selecting the rows of the identity matrix, it is clear that every row of both ${\bf C}_t \otimes {\bf C}_t$ and ${\bf T}$ only contains a single one and zeros elsewhere. 
This fact guarantees that each row of ${\bf R}_{c_t}$ has only a single one and thus, in order to ensure the full column rank condition of ${\bf R}_{c_t}$, we need to ensure that each column of it has at least a single one. This problem actually has been encountered and solved in~\cite{TSP12}  
where the solution is to construct ${\bf C}_t$ by selecting the rows of ${\bf I}_{N_t}$ based on the so-called minimal length-$(N_t-1)$ sparse ruler problem. In practice, this results in a multi-coset sampling procedure called the minimal sparse ruler sampling~\cite{TSP12}. 

Next, we examine the choice of 
${\bf B}$, which boils down to the selection of the activated antennas in the ULA and the investigated angles $\{{\theta}_q\}_{q=1}^{Q}$. Let us write ${\bf B}^*\odot{\bf B}$ in terms of $\{{\bf b}({\theta}_q)\}_{q=1}^Q$ as
\vspace{-1mm}
\begin{equation}
{\bf B}^*\odot{\bf B}=\left[{\bf b}^*({\theta}_1)\otimes{\bf b}({\theta}_1),\dots,{\bf b}^*({\theta}_Q)\otimes{\bf b}({\theta}_Q)\right]
\vspace{-1mm}
\label{eq:Btilde_as_btilde}
\end{equation}
and ${\bf b}({\theta}_q)$ in terms of $a(\theta_q)$ as 
\vspace{-1mm}
\begin{equation}
{\bf b}({\theta}_q)=\left[a({\theta}_q)^{{d}_1}, a({\theta}_q)^{{d}_2}, \dots, a({\theta}_q)^{{d}_{M_s}}\right]^T
\vspace{-1mm}
\label{eq:btilde}
\end{equation}
where ${d}_i$ is the distance in wavelengths between the $i$-th {\it active} antenna and the reference antenna of the ULA defined in Section~\ref{preliminary}. 
It is clear from~\eqref{eq:Btilde_as_btilde} and~\eqref{eq:btilde} that the $q$-th column of ${\bf B}^*\odot{\bf B}$ contains the elements $\text{exp}\left(j ({d}_i-{d}_j) 2\pi\text{sin}({\theta}_q)\right)$, for $i,j=1,2,\dots,M_s$. While our task to find general design conditions to guarantee the full column rank of ${\bf B}^*\odot{\bf B}$ is not trivial, the following theorem suggests one possible way to achieve a full column rank ${\bf B}^*\odot{\bf B}$.
\vspace*{1mm}
\newline
\hspace*{1mm}\textbf{Theorem 1}: The matrix ${\bf B}^*\odot{\bf B}$ has full column rank if: 
1) There exist $Q$ distinct values of $\theta_q$ satisfying $-\frac{\pi}{2}<\{\theta_q\}_{q=1}^Q\leq\frac{\pi}{2}$, and 
2) There exists an integer $N_v\geq Q$ such that $\{{d}_i-{d}_j\}_{i,j=1}^{M_s}$ contains an arithmetic sequence of $N_v$ terms having a difference of $d\leq 0.5$ between each two consecutive terms.
\vspace*{1mm}
\newline
The proof of Theorem~1 can be found in Appendix~\ref{Appendix_fullrank_B_kr_B}. The second condition indicates that there exist $N_v$ distinct rows from ${\bf B}^*\odot{\bf B}$ that form the array response matrix of a virtual ULA with $N_v$ antennas, which can only be achieved for $N_v \leq 2N_s-1$. This second condition also implies that we have more antennas in this virtual ULA than investigated angles. Some possible ways to satisfy Theorem 1 is to select the $M_s$ active antennas from the $N_s$ antennas in the ULA based on the MRA discussed in~\cite{Moffet} (which also obeys the minimal sparse ruler problem~\cite{Siavash}), the two-level nested array~\cite{Nested}, or the coprime array~\cite{Coprime}. For the MRA and the two-level nested array, Theorem 1 can be satisfied even for $N_v=2N_s-1$. Note that although the $Q$ different values of ${\theta}_q$ can be chosen in an arbitrary fashion, they should not be too close to each other, since otherwise the resulting ${\bf B}^*\odot{\bf B}$ might be ill-conditioned. Theorem $1$ also implies that the maximum number of detectable sources is upper bounded by $K\leq 2N_s-1$ since we cannot detect more than $Q$ sources.
Apart from satisfying Theorem 1, another way to achieve a full column rank ${\bf B}^*\odot{\bf B}$ is suggested by Theorem 2.
\vspace*{1mm}
\newline
\hspace*{1mm}\textbf{Theorem 2}: The matrix ${\bf B}^*\odot{\bf B}$ has full column rank if:\newline
1) $\{({d}_i-{d}_j) \, \mathrm{mod} \, \frac{Q}{2}\}_{i,j=1}^{M_s}$ has at least $Q$ different values and 2) the grid of investigated angles $\{{\theta}_q\}_{q=1}^Q$ is designed based on the inverse sinusoidal angular grid where 
\vspace{-1mm}
\begin{equation}
\small
{\theta}_q=\text{sin}^{-1}\left(\frac{2}{Q}\left(q-1-\left\lceil\frac{Q-1}{2}\right\rceil\right)\right),
\label{eq:inverse_sin_grid}
\end{equation} 
The proof for this theorem can be found in Appendix~\ref{Appendix_fullrank_B_kr_B_Two}. Note that the 
first condition from Theorem~2 is less strict than the second condition from Theorem~1.
A good option is to use a configuration satisfying Theorem 1 with $N_v=2N_s-1$ and $d=0.5$, and to use~\eqref{eq:inverse_sin_grid} 
with $Q=2N_s-1$. This will not only ensure that the resulting $M_s^2 \times (2N_s-1)$ matrix ${\bf B}^*\odot{\bf B}$ has full column rank but also that there exists a $(2N_s-1) \times (2N_s-1)$ submatrix from ${\bf B}^*\odot{\bf B}$ that forms a row-permuted version of the $(2N_s-1) \times (2N_s-1)$ inverse DFT matrix, meaning that ${\bf B}^*\odot{\bf B}$ is well-conditioned.

\vspace{-1mm}
\section{Numerical Study}
\vspace{-1mm}
In this section, we examine the proposed approach with some numerical study. We consider a ULA having $N_s=36$ antennas as the underlying array and construct an MRA of active antennas by selecting the antenna indices based on the minimal length-$35$ sparse ruler problem discussed in~\cite{TSP12,Siavash}. This leads to $M_s=10$ activated antennas with $\{{d}_j\}_{j=1}^{10}=$ $\{0,d,4d,10d,16d,22d,28d,30d,33d,35d\}$ where $d$ is set to $d=0.5$. 
The set of investigated angles $\{{\theta}_q\}_{q=1}^Q$ is set according to~\eqref{eq:inverse_sin_grid} with $Q=2N_s-1=71$.
In the receiver branch corresponding to each active antenna, the time-domain compression rate of $M_t/N_t=0.4048$ 
is obtained by setting $N_t=84$ and $M_t=34$. 
We construct the $34 \times 84$ selection matrix ${\bf C}_t$ by first solving the minimal length-$83$ sparse ruler problem which gives the indices of the $16$ rows of ${\bf I}_{N_t}$ that have to be selected. The selection of these $16$ rows will ensure that the resulting matrix ${\bf R}_{c_t}$ in~\eqref{eq:vec_R_yi_yj_prime_as_T} 
has at least a single one in each column. The additional $18$ rows of ${\bf C}_t$ are then 
randomly selected from the remaining rows of ${\bf I}_{N_t}$ that have not been selected.
We simulate the case when we have more sources than active antennas by generating $K=12$ uncorrelated sources 
having DOAs 
with 9 degrees of separation, i.e., the set of DOAs is given by 
$\{-54^0,-45^0,\dots,45^0\}$. The sources 
produce complex baseband signals whose frequency bands are given in Table~\ref{tab:mytabel} and which are 
generated by passing circular complex zero-mean Gaussian i.i.d. noise with variance $\sigma^2=5$ into a digital filter of length $N_t=84$ with the unit-gain passband of the filter for each source set according to Table~\ref{tab:mytabel}. This will ensure that the true auto-correlation sequence 
for each source 
is limited to $-N_t+1\leq m \leq N_t-1$. 
We assume a spatially and temporally white noise 
with variance $\sigma_n^2=5$ and set 
the number of measurement matrices ${\bf Z}[n]$ to $N_n=5951$.
\vspace{-2mm}
\begin{table}[h]
\small
	\caption{The frequency band occupied by the sources}
	\centering
		\begin{tabular}{| c | c | c |}
			\hline
    Source & Actual DOA & Occupied frequency band \\ \hline
    $1$ & $-54^0$ & $[-0.275\pi,-0.2\pi]$ \\ \hline
    $2$ & $-45^0$ & $[-0.8\pi,-0.725\pi]$ \\ \hline
    $3$ & $-36^0$ & $[-0.35\pi,-0.275\pi]$ \\ \hline
    $4$ & $-27^0$ & $[0.35\pi,0.425\pi]$ \\ \hline
    $5$ & $-18^0$ & $[0.875\pi,0.95\pi]$ \\ \hline
    $6$ & $-9^0$ & $[0.05\pi,0.125\pi]$ \\ \hline
    $7$ & $0^0$ & $[-0.95\pi,-0.875\pi]$ \\ \hline
    $8$ & $9^0$ & $[-0.65\pi,-0.575\pi]$ \\ \hline
    $9$ & $18^0$ & $[-0.425\pi,-0.35\pi]$ \\ \hline
    $10$ & $27^0$ & $[0.575\pi,0.65\pi]$ \\ \hline
    $11$ & $36^0$ & $[0.125\pi,0.2\pi]$ \\ \hline
    $12$ & $45^0$ & $[0.5\pi,0.575\pi]$ \\ \hline
    \end{tabular}
\label{tab:mytabel}
\end{table}

\begin{figure}[h]
    \centering
        \includegraphics[width=0.5\textwidth]{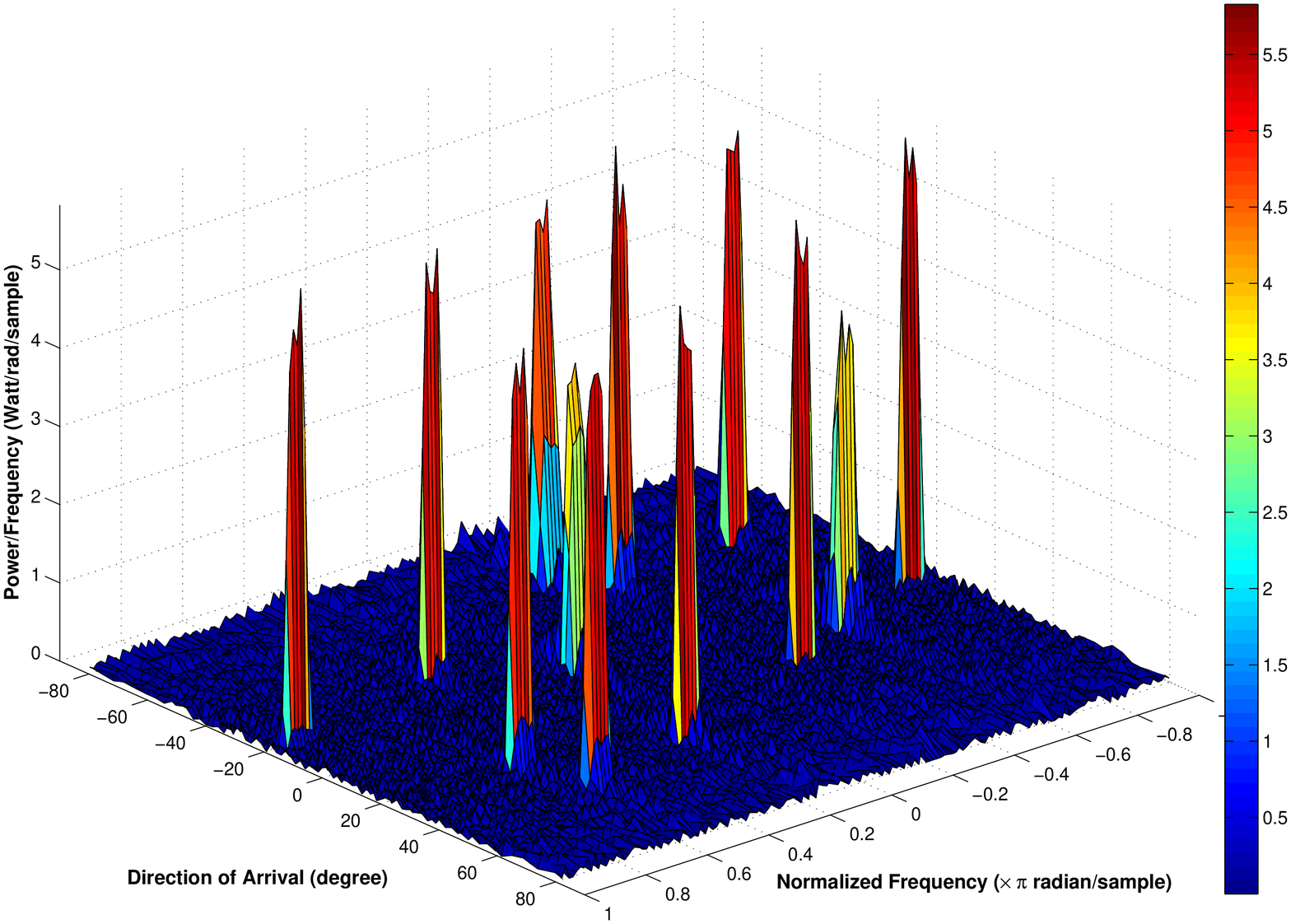}
        \vspace{-1cm}
        \caption{\small The power spectrum estimate (in watt/radian/sample) as a function of frequency (radian/sample) and angle (degree). 
        }
    \label{fig:Normal_plots_2DPSD}
    \centering
        \includegraphics[width=0.5\textwidth]{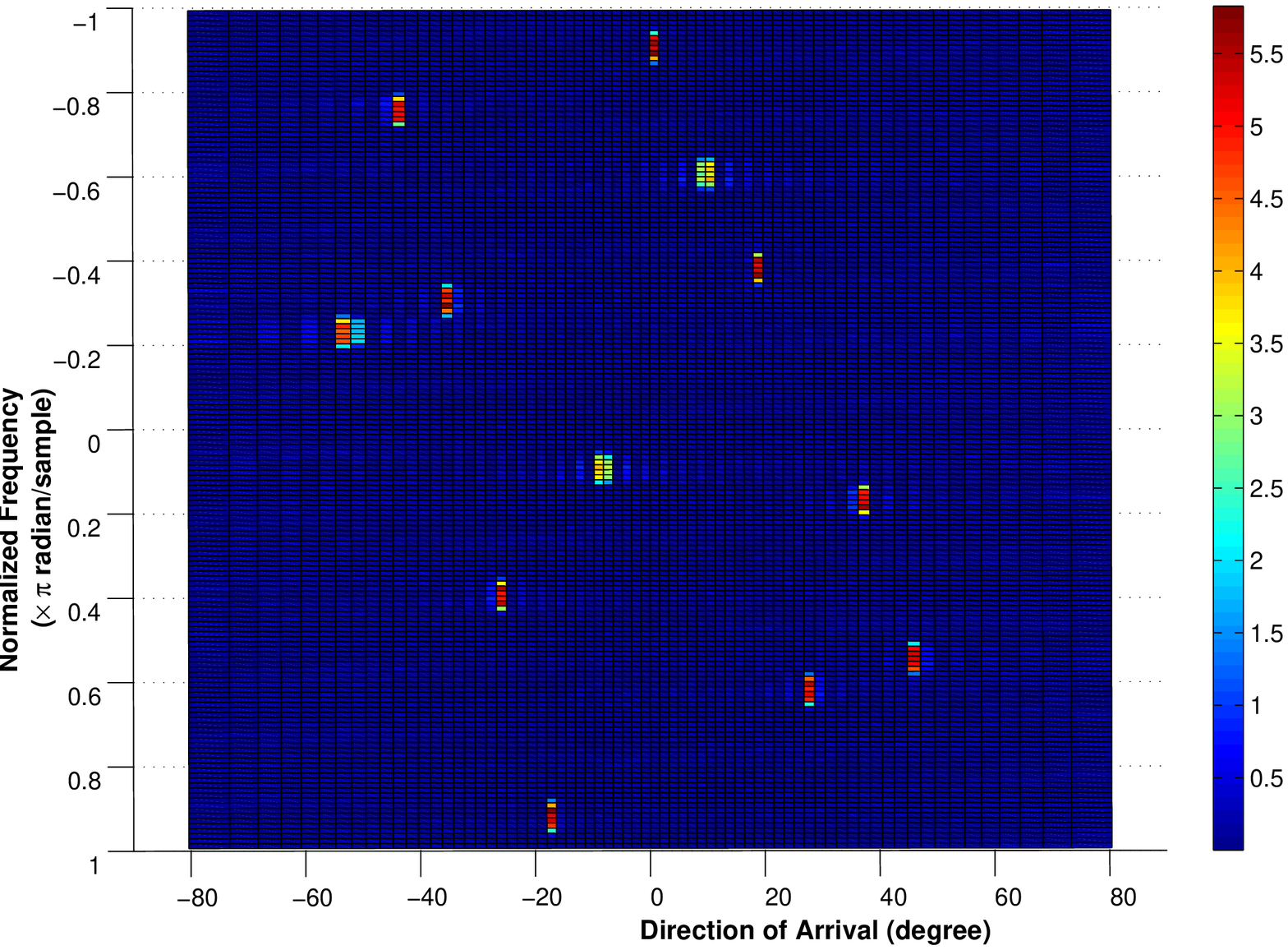}
        \vspace{-0.7cm}
        \caption{\small The top view of Fig.~\ref{fig:Normal_plots_2DPSD}.} 
    \label{fig:Top_plots_2DPSD}
\end{figure}

Fig.~\ref{fig:Normal_plots_2DPSD} illustrates the estimate of the power spectrum as a function of the frequency and the investigated angles. It is clear that the $12$ uncorrelated sources can generally be detected. We can find the DOA estimates by locating the peak of this spectrum though the actual DOAs might not fall on top of the defined investigated angles. For a given DOA estimate, we can locate the active frequency band of the corresponding source together with the value of the power spectrum estimate. The top view of Fig.~\ref{fig:Normal_plots_2DPSD}, which is provided by Fig.~\ref{fig:Top_plots_2DPSD}, gives a much clearer picture of the quality of the estimate. We can easily compare this figure with the data provided in Table~\ref{tab:mytabel}. Observe that the estimate of the DOA, the power spectrum, as well as the active frequency band of the sources is quite satisfactory except for the sources with DOAs 
$-9^0$ and 
$9^0$. For these two sources, it is apparent from Fig.~\ref{fig:Top_plots_2DPSD} that the impact of the grid mismatch effect is quite significant and their power spectrum estimates seem to have been distributed among the two nearest grid points. Note that this 2D power spectrum estimate can be produced without applying any sparsity contraint on the true power spectrum, but can of course be improved if such a constraint is used. 
%

\appendix
\vspace{-1mm}
\section{Proof of Theorem 1}
\vspace{-1mm}
\label{Appendix_fullrank_B_kr_B}
The second requirement of Theorem~1 implies that there exists a $Q \times Q$ matrix $\acute{\bf B}=[\acute{\bf b}({\theta}_1),\acute{\bf b}({\theta}_2),\dots,\acute{\bf b}({\theta}_Q)]$, which is a submatrix of ${\bf B}^*\odot{\bf B}$ in~\eqref{eq:Btilde_as_btilde}, that forms the array response matrix of a virtual ULA of $Q$ antennas with $\acute{\bf b}({\theta}_q)$ given by
$\acute{\bf b}({\theta}_q)=[a({\theta}_q)^{\bar{d}},a({\theta}_q)^{\bar{d}+d},\dots,a({\theta}_q)^{\bar{d}+(Q-1)d}]^T$, where ${\bar{d}}$ gives the distance between the first antenna in the virtual ULA and the reference antenna in the underlying ULA in Section~\ref{preliminary}. 
Hence, it is clear that $\acute{\bf B}$ is a column-wise Vandermonde matrix. 
From the well-known properties of a column-wise Vandermonde matrix, $\acute{\bf B}$ has full column rank due to the first requirement of Theorem~1
and since $d\leq 0.5$. It is then trivial to show that ${\bf B}^*\odot{\bf B}$ also has full column rank.
\vspace{-1mm}
\section{Proof of Theorem 2}
\vspace{-1mm}
\label{Appendix_fullrank_B_kr_B_Two}
Based on~\eqref{eq:Btilde_as_btilde} and~\eqref{eq:btilde} and the fact that the inverse sinusoidal angular grid in~\eqref{eq:inverse_sin_grid} is used, we can write ${\bf B}^*\odot{\bf B}$ in terms of its row vectors, i.e., ${\bf B}^*\odot{\bf B} = [{\boldsymbol \beta}({d}_1-{d}_1),{\boldsymbol \beta}({d}_2-{d}_1),\dots,{\boldsymbol \beta}({d}_{M_s}-{d}_{M_s})]^T$, with ${\boldsymbol \beta}({d}_{i}-{d}_{j})$ given by
${\boldsymbol \beta}({d}_{i}-{d}_{j})=[e^{j\frac{4\pi}{Q}({d}_i-{d}_j)(-\left\lceil\frac{Q-1}{2}\right\rceil)},\dots,e^{j\frac{4\pi}{Q}({d}_i-{d}_j)(-1)},$ 
$1,e^{j\frac{4\pi}{Q}({d}_i-{d}_j)},
\dots,e^{j\frac{4\pi}{Q}({d}_i-{d}_j)(Q-1-\left\lceil\frac{Q-1}{2}\right\rceil)}]^T$.
Observe that ${\bf B}^*\odot{\bf B}$ is a row-wise Vandermonde matrix since the elements of ${\boldsymbol \beta}({d}_{i}-{d}_{j})$ are ordered according to geometric progression. 
In order to ensure that ${\bf B}^*\odot{\bf B}$ has full column rank, we need $Q$ distinct values of $\frac{4\pi}{Q}({d}_i-{d}_j)$ modulo $2\pi$
which 
is guaranteed by the first requirement of Theorem 2. 

\vspace{-5mm}
\vspace{-0.5mm}
\end{document}